\documentclass[11pt]{amsart}
\usepackage{eucal}

\newtheorem{thm}{Theorem}

\newtheorem{prop}[thm]{Proposition}

\newcommand{\R}{\mathbb{R}}
\newcommand{\E}{\mathbb{E}}
\newcommand{\Prob}{\mathbb{P}}

\newcommand{\Z}{\mathbb{Z}}
\newcommand{\C}{\mathbb{C}}
\newcommand{\I}{{\bf 1}}

\newcommand{\inprod}[2]{\left\langle #1, #2 \right\rangle}

\allowdisplaybreaks[4]

\title{On the spectral norm of a random Toeplitz matrix}

\author[M.\ Meckes]{Mark W.\ Meckes}

\email{mwmeckes@math.cornell.edu}

\address{Department of Mathematics, Cornell University, Ithaca,
New York 14853, U.S.A.}

\begin{document}

\begin{abstract}
Suppose that $T_n$ is a Toeplitz matrix whose entries come from a
sequence of independent but not necessarily identically distributed
random variables with mean zero. Under some additional tail
conditions, we show that the spectral norm of $T_n$ is of the order
$\sqrt{n \log n}$. The same result holds for random Hankel matrices
as well as other variants of random Toeplitz matrices which have
been studied in the literature.
\end{abstract}

\maketitle


\section{Introduction and results}\label{S:intro}

Let $X_0,X_1,X_2,\dotsc$ be a family of independent random
variables. For $n\ge 2$, $T_n$ denotes the $n\times n$ random
symmetric Toeplitz matrix $T_n = \big[ X_{|j-k|}\big]_{1\le j, k\le
n}$,
\[
T_n = \begin{bmatrix}
X_0 & X_1 & X_2 & \cdots & X_{n-2} & X_{n-1} \\
X_1 & X_0 & X_1 & & & X_{n-2} \\
X_2 & X_1 & X_0 & & & \vdots \\
\vdots & & & \ddots & & \vdots\\
X_{n-2} & & & & X_0 & X_1 \\
X_{n-1} & X_{n-2} & \hdotsfor{2} & X_1 & X_0
\end{bmatrix}.
\]
In \cite{Bai}, Bai asked whether the spectral measure of
$n^{-1/2} T_n$ approaches a deterministic limit measure $\mu$
as $n\to \infty$. Bryc, Dembo, and Jiang \cite{BDJ} and Hammond and
Miller \cite{HM} independently proved that this is so when the $X_j$
are identically distributed with variance $1$, and that with these
assumptions $\mu$ does not depend on the distribution of the
$X_j$. The measure $\mu$ does not appear to be a previously studied
probability measure, and is described via rather complicated
expressions for its moments.

This limiting spectral measure $\mu$ has unbounded support, which
raises the question of the asymptotic behavior of the spectral norm
$\| T_n \|$, i.e., the maximum absolute value of an eigenvalue of
$T_n$. (This problem is explicitly raised in \cite[Remark 1.3]{BDJ}.)
This paper shows, under slightly different assumptions from
\cite{BDJ,HM}, that $\| T_n \|$ is of the order $\sqrt{n\log n}$. Here
the $X_j$ need not be identically distributed, but satisfy stronger
moment or tail conditions than in \cite{BDJ,HM}. The spectral norm is
also of the same order for other related random matrix ensembles,
including random Hankel matrices. In the case of Hankel matrices,
Theorems \ref{T:upper} and \ref{T:lower} below generalize in a
different direction a special case of a result of Masri and Tonge
\cite{HM} on multilinear Hankel forms with $\pm 1$ Bernoulli entries.

\medskip

A random variable $X$ will be called \emph{subgaussian} if
\begin{equation}\label{E:subgaussian}
\Prob \big[|X| \ge t\big] \le 2 e^{-at^2} \quad \forall t>0
\end{equation}
for some constant $a>0$. A family of random variables is
\emph{uniformly subgaussian} if each satisfies \eqref{E:subgaussian}
for the same constant $a$.

\begin{thm}\label{T:upper}
Suppose $X_0,X_1,X_2,\dotsc$ are independent, uniformly subgaussian
random variables with $\E X_j = 0$ for all $j$. Then
\[
\E \| T_n \| \le c_1 \sqrt{n\log n},
\]
where $c_1 > 0$ depends only on the constant $a$ in the subgaussian estimate
\eqref{E:subgaussian} for the $X_j$.
\end{thm}

Simple scaling considerations show that one can take $c_1 = C
a^{-1/2}$ for some absolute constant $C>0$. In principle an explicit
value for $C$ can be extracted from the proof of Theorem
\ref{T:upper}.  No attempt has been made to do so, since the
techniques used in this paper are suited for determining rough orders
of growth, and not precise constants. Similar remarks apply to the
constants which appear in the statements of Theorems \ref{T:limsup}
and \ref{T:lower} below.

\medskip

By strengthening the subgaussian assumption, the statement of Theorem
\ref{T:upper} can be improved from a bound on expectations to an
almost sure asymptotic bound. Recall that a real-valued random
variable $X$ (or more properly, its distribution) is said to satisfy a
\emph{logarithmic Sobolev inequality} with constant $A$ if
\[
\E \big[f^2(X) \log f^2(X)\big] \le 2A \ \E \big[ f'(X)^2\big]
\]
for every smooth $f:\R\to \R$ such that $\E f^2(X)=1$. Standard normal
random variables satisfy a logarithmic Sobolev inequality with
constant $1$. Furthermore, it is well known that independent random
variables with bounded logarithmic Sobolev constants are uniformly
subgaussian and possess the same concentration properties as
independent normal random variables (see \cite{Ledoux1} or
\cite[Chapter 5]{Ledoux2}).

\begin{thm}\label{T:limsup}
Suppose $X_0,X_1,X_2,\dotsc$ are independent, $\E X_j=0$ for all $j$,
and for some constant $A$, either:
\begin{enumerate}
\item \label{I:bounded} for all $j$, $|X_j|\le A$ almost surely; or
\item \label{I:LSI} for all $j$, $X_j$ satisfies a logarithmic Sobolev
inequality with constant $A$.
\end{enumerate}
Then 
\[
\limsup_{n\to\infty} \frac{\| T_n \|}{\sqrt{n\log n}} \le c_2
\]
almost surely, where $c_2 > 0$ depends only on $A$.
\end{thm}

We remark that according to the definition used here, $T_n$ is a
submatrix of $T_{n+1}$, but this is only a matter of convenience in
notation. Theorem \ref{T:limsup} remains true regardless of the
dependence among the random matrices $T_n$ for different values of
$n$.

It seems unlikely that the stronger hypotheses of Theorem
\ref{T:limsup} are necessary. In fact a weaker version can be proved
under the hypotheses of Theorem \ref{T:upper} alone; see the remarks
following the proof of Theorem \ref{T:limsup} in Section
\ref{S:proofs}.

\medskip

When the $X_j$ have variance $1$, the upper bound $\sqrt{n\log n}$ of
Theorems \ref{T:upper} and \ref{T:limsup} is of the correct order. In
fact the matching lower bound holds under less restrictive tail
assumptions, as the next result shows.

\begin{thm}\label{T:lower}
Suppose $X_0,X_1,X_2,\dotsc$ are independent and for some constant $B$,
each $X_j$ satisfies
\[
\E X_j = 0, \quad \E X_j^2 = 1, \quad \E |X_j| \ge B.
\]
Then
\[
\E \| T_n \| \ge c_3 \sqrt{n\log n},
\]
where $c_3 > 0$ depends only on $B$.
\end{thm}

In the case that $\E X_j^2 = 1$ and $\E |X_j|^3 < \infty$, it is a
consequence of H\"older's inequality that $\E |X_j| \ge (\E
|X_j|^3)^{-1}$.  Thus the lower bound on first absolute moments
assumed in Theorem \ref{T:lower} is weaker than an upper bound on
absolute third moments, and is in particular satisfied for uniformly
subgaussian random variables.

\medskip

Section \ref{S:proofs} below contains the proofs of Theorems
\ref{T:upper}--\ref{T:lower}.As mentioned above, Theorems
\ref{T:upper}--\ref{T:lower} also hold for other ensembles of random
Toeplitz matrices, as well as for random Hankel matrices.  Section
\ref{S:remarks} discusses these extensions of the theorems and makes
some additional remarks.

\medskip

{\em Acknowledgement.} The author thanks A.\ Dembo for pointing out
the problem considered in this paper.

\bigskip
\section{Proofs}\label{S:proofs}

The proof of Theorem \ref{T:upper} is based on Dudley's entropy bound
\cite{Dudley} for the supremum of a subgaussian random process. Given
a random process $\{Y_x : x\in M\}$, a pseudometric on $M$ may be
defined by
\[
d(x,y) = \sqrt{\E |Y_x - Y_y|^2}.
\]
The process $\{Y_x : x\in M\}$ is called \emph{subgaussian} if
\begin{equation}\label{E:process}
\forall x,y\in M,\ \forall t>0, \quad 
\Prob \big[|Y_x - Y_y| \ge t \big] \le 2 \exp\left[ -
  \frac{b\ t^2}{d(x,y)^2} \right]
\end{equation}
for some constant $b>0$. For $\varepsilon > 0$, the $\varepsilon$-covering
number of $(M,d)$, $N(M,d,\varepsilon)$, is the smallest cardinality of a
subset $\mathcal{N}\subset M$ such that
\[
\forall x\in M \ \exists y\in \mathcal{N} : \ d(x,y)\le \varepsilon.
\]
Dudley's entropy bound is the following (see \cite[Proposition
2.1]{Talagrand2} for the version given here).

\begin{prop}\label{T:Dudley}
Let $\{Y_x : x\in M\}$ be a subgaussian random process with $\E Y_x =
0$ for every $x\in M$. Then
\[
\E \sup_{x\in M} |Y_x| \le K \int_0^\infty \sqrt{\log N(M,d,\varepsilon)}
  \ d \varepsilon,
\]
where $K>0$ depends only on the constant $b$ in the subgaussian estimate
\eqref{E:process} for the process.
\end{prop}

\medskip

We will also need the following version of the classical Azuma-Hoeffding
inequality. This can be proved by a standard Laplace transform argument;
see e.g.\ \cite[Fact 2.1]{LPRT}.

\begin{prop}\label{T:AH}
Let $X_1,\dotsc,X_n$ be independent, symmetric, uniformly subgaussian
random variables. Then for any $a_1,\dotsc,a_n \in \R$ and $t>0$,
\[
\Prob \Biggl[\biggl|\sum_{j=1}^n a_j X_j \biggr| \ge t \Biggr]
  \le 2 \exp \left[-\frac{b\ t^2}{\sum_{j=1}^n a_j^2}\right],
\]
where $b>0$ depends only on the constant $a$ in the subgaussian estimate
\eqref{E:subgaussian} for the $X_j$.
\end{prop}

\medskip

\begin{proof}[Proof of Theorem \ref{T:upper}]
We first reduce to the case in which each $X_j$ is symmetric. Let
$T_n'$ be an independent copy of $T_n$. Since $\E T_n = 0$, by
Jensen's inequality,
\[
\E \| T_n \| \le \E \big[ \E \big[ \|T_n-T_n'\| \big| T_n \big]\big]
  = \E \|T_n-T_n'\|.
\]
The random Toeplitz matrix $(T_n-T_n')$ has entries $(X_j-X_j')$ which
are independent, symmetric, uniformly subgaussian random variables
(with a possibly smaller constant $a$ in the subgaussian estimate).
Thus we may assume without loss of generality that the $X_j$ are
symmetric random variables.

\medskip

We next bound $\|T_n\|$ by the supremum of a subgaussian random
process. A basic feature of the theory of Toeplitz matrices is
their relationship to multiplication operators (cf.\ \cite[Chapter
1]{BS}). Specifically, the finite Toeplitz matrix $T_n$ is an $n\times
n$ submatrix of the infinite Laurent matrix
\[
L_n = \big[ X_{|j-k|} \I_{|j-k| \le n-1} \big]_{j,k \in \Z}.
\]
Consider $L_n$ as an operator on $\ell^2(\Z)$ in the canonical way,
and let $\psi:\ell^2(\Z) \to L^2[0,1]$ denote the usual trigonometric
isometry $\psi(e_j)(x) = e^{2\pi ijx}$.  Then $\psi L_n
\psi^{-1}:L^2\to L^2$ is the multiplication operator corresponding to
the $L^\infty$ function
\[
f(x) = \sum_{j=-(n-1)}^{n-1} X_{|j|} e^{2\pi i jx}
     = X_0 + 2 \sum_{j=1}^{n-1} \cos(2\pi j x) X_j.
\]
Therefore
\begin{equation}\label{E:norms}
\|T_n \| \le \|L_n\| = \| f \|_\infty = \sup_{0\le x \le 1} |Y_x|,
\end{equation}
where
\[
Y_x = X_0 + 2 \sum_{j=1}^{n-1} \cos(2\pi j x) X_j.
\]
By Proposition \ref{T:AH}, the random process $\{Y_x : x\in[0,1]\}$
becomes subgaussian if $M=[0,1]$ is equipped with the pseudometric
\[
d(x,y) = \sqrt{\sum_{j=1}^{n-1} \big[ \cos(2\pi j x) 
  - \cos(2\pi jy)\big]^2}.
\]

\medskip

Finally, we bound $N([0,1],d,\varepsilon)$ in order to apply Proposition
\ref{T:Dudley}. Since $|\cos t| \le 1$ always, it follows that $d(x,y)
< 2\sqrt{n}$ and therefore $N([0,1],d,\varepsilon) = 1$ if $\varepsilon
> 2\sqrt{n}$.  Next, since $|\cos s - \cos t| \le |s-t|$, 
\[
d(x,y) \le 2\pi |x-y| \sqrt{\sum_{j=1}^{n-1} j^2} < 4 n^{3/2} |x-y|,
\]
which implies that
\[
N\bigl([0,1],d,\varepsilon\bigr) \le N\left([0,1],|\cdot |, 
  \frac{\varepsilon}{4n^{3/2}}\right) \le \frac{4n^{3/2}}{\varepsilon}.
\]
By \eqref{E:norms}, Proposition \ref{T:Dudley}, and the substitution
$\varepsilon = 4n^{3/2}e^{-t^2}$,
\begin{equation}\label{E:upper-Dudley}
\E \|T_n\| \le K \int_0^{2\sqrt{n}} \sqrt{\log\left(
  \frac{4n^{3/2}}{\varepsilon}\right)}\ d\varepsilon
  = 2\sqrt{2} n^{3/2} K \int_{\sqrt{2\log 2n}}^\infty t^2 e^{-t^2/2} \ dt.
\end{equation}
Integration by parts and the classical estimate
$\frac{1}{\sqrt{2\pi}}\int_s^\infty e^{-t^2/2}\ dt \le e^{-s^2/2}$ for
$s>0$ yield
\[
\int_s^\infty t^2 e^{-t^2/2}\ dt \le \big(s+\sqrt{2\pi}\big) e^{-s^2/2}.
\]
Combining the case $s=\sqrt{2\log 2n}$ of this estimate with
\eqref{E:upper-Dudley} completes the proof.
\end{proof}

\medskip

The proof of Theorem \ref{T:limsup} is based on rather classical
measure concentration arguments commonly applied to probability in
Banach spaces.

\begin{proof}[Proof of Theorem \ref{T:limsup}]
Denote by $M_0$ the $n\times n$ identity matrix, and for
$m=1,\dotsc,n-1$ let $M_m = \bigl[\I_{|j-k|=m}\bigr]_{1\le j,k\le
n}$. Then $T_n$ can be written as the sum
\[
T_n = \sum_{j=0}^{n-1} X_j M_j
\]
of independent random vectors in the finite-dimensional Banach space
$\mathcal{M}_n$ equipped with the spectral norm. Observe that $\|M_j\|
\le 2$ for every $j$.

Under the assumption (\ref{I:bounded}), up to the precise values of
constants the estimate
\[
\Prob \bigl[ \| T_n \| \ge \E \| T_n \| + t \bigr]
  \le e^{-t^2 / 32A^2n} \quad \forall t>0
\]
follows from any of several standard approaches to concentration of
measure (cf.\ Corollary 1.17, Corollary 4.5, or Theorem 7.3 of
\cite{Ledoux2}; the precise statement can be proved from Corollary
1.17). Combining this with Theorem 1 yields
\[
\Prob\Bigl[\|T_n \| \ge (c_1 + 8A) \sqrt{n\log n}\Bigr] \le \frac{1}{n^2},
\]
which completes the proof via the Borel-Cantelli lemma.

The proof under the assumption (\ref{I:LSI}) is similar. By the triangle
inequality and the Cauchy-Schwarz inequality,
\[
\|T_n\| \le 2 \sqrt{n \sum_{j=0}^{n-1} X_j^2},
\]
so that the map $(X_0,\dotsc,X_{n-1}) \mapsto \|T_n\|$ has Lipschitz
constant bounded by $2\sqrt{n}$. By the well-known tensorization and
measure concentration properties of logarithmic Sobolev inequalities
(cf.\ \cite[Sections 2.1--2.3]{Ledoux1} or \cite[Sections
5.1--5.2]{Ledoux2}),
\[
\Prob \bigl[ \| T_n \| \ge \E \| T_n \| + t \big]
  \le e^{-t^2 / 4An} \quad \forall t>0.
\]
The proof is completed in the same way as before (with a different
dependence of $c_2$ on $A$).
\end{proof}

As remarked above, a weaker version of Theorem \ref{T:limsup} may be
proved under the assumptions of Theorem \ref{T:upper} alone. From the
proof of Proposition \ref{T:Dudley} in \cite{Talagrand2} one can
extract the following tail inequality under the assumptions of
Proposition \ref{T:Dudley}:
\begin{equation}\label{E:Dudley-tail}
\Prob \biggl[\sup_{x\in M} |Y_x| \ge t\biggr] \le 2 e^{-c t^2/\alpha^2}
\quad \forall t>0, \quad \mbox{where} \quad
\alpha=\int_0^\infty \sqrt{\log N(M,d,\varepsilon)} \ d \varepsilon.
\end{equation}
The explicit statement here is adapted from lecture notes of Rudelson
\cite{Rudelson}. Using the estimates derived in the proof of Theorem
\ref{T:upper} and applying the Borel-Cantelli lemma as above, one
directly obtains
\begin{equation}\label{E:weak-limsup}
\limsup_{n\to\infty} \frac{\| T_n \|}{\sqrt{n}\log n} \le c_4
\quad \mbox{almost surely}
\end{equation}
under the assumptions that the $X_j$ are symmetric and uniformly
subgaussian. The general (nonsymmetric but mean $0$) case can be
deduced from the argument for the symmetric case.  Let $T_n'$ be an
independent copy of $T_n$. By independence, the triangle inequality,
and the tail estimate which follows from \eqref{E:Dudley-tail},
\[
\Prob\big[ \|T_n'\| \le s \big] \Prob\big[\|T_n\| \ge s+t\big]
\le \Prob \big[ \|T_n - T_n'\| \ge t\big]
\le 2 e^{-ct^2/n\log n}
\]
for some constant $c$ which depends on the subgaussian estimate for
the $X_j$. By Theorem \ref{T:upper} and Chebyshev's
inequality,
\[
\Prob\big[\|T_n'\| \le s \big] \ge 1-\frac{1}{s}c_1 \sqrt{n\log n}.
\]
Picking $s=2 c_1 \sqrt{n\log n}$ and $t=\sqrt{\frac{2n}{c}}\log n$
yields
\[
\Prob[ \big[\| T_n \| \ge c_4 \sqrt{n}\log n \big] \le \frac{4}{n^2}
\]
for some constant $c_4$, and \eqref{E:weak-limsup} then follows from
the Borel-Cantelli lemma.

\medskip

The proof of Theorem \ref{T:lower} amounts to an adaptation of the
proof of the lower bound in \cite{MT}, with much of the proof
abstracted into a general lower bound for the suprema of certain
random processes due to Kashin and Tzafriri \cite{KT1,KT2}.  The
following is a special case of the result of \cite{KT2}.

\begin{prop}\label{T:KT}
Let $\varphi_j:[0,1]\to \R$, $j=0,\dotsc,n-1$ be a family of functions
which are orthonormal in $L^2[0,1]$ and satisfy $\| \varphi_j \|_{L^3[0,1]}
\le A$ for every $j$, and let $X_0,\dotsc,X_{n-1}$ be independent random
variables such that for every $j$,
\[
\E X_j = 0, \quad \E X_j^2 = 1, \quad \E |X_j| \ge B.
\]
Then for any $a_0,\dotsc,a_{n-1}\in \R$,
\[
\E \left[ \sup_{0\le x \le 1}
  \Biggl|\sum_{j=0}^{n-1} a_j X_j \varphi_j(x) \Biggr| \right]
\ge K\ \|a\|_2 \sqrt{\log \frac{\|a\|_2}{\|a\|_4}},
\]
where $\|a\|_p = \big(\sum_{j=0}^{n-1} |a_j|^p\big)^{1/p}$ and $K>0$
depends only on $A$ and $B$.
\end{prop}

\medskip

\begin{proof}[Proof of Theorem \ref{T:lower}]
First make the estimate
\[
\|T_n\| = \sup_{v\in\C^n\setminus \{0\} } 
  \frac{|\inprod{T_n v}{v}|}{\inprod{v}{v}}
  \ge \sup_{0\le x \le 1}\frac{1}{n} \big\vert\inprod{T_n v_x}{v_x}\big\vert,
\]
where $v_x\in \C^n$ is defined by $(v_x)_j = e^{2\pi i j x}$ for
$j=1,\dotsc,n$ and $\inprod{\cdot}{\cdot}$ is the standard inner
product on $\C^n$. Therefore
\begin{align*}
\|T_n\| & \ge \frac{1}{n} \sup_{0\le x\le 1} \Biggl| 
  \sum_{j,k=1}^n X_{|j-k|} e^{2\pi i (j-k)x} \Biggr|\\
&= \frac{1}{n} \sup_{0\le x\le 1}\Biggl| \sum_{j=-(n-1)}^{n-1} (n-|j|) X_{|j|}
  e^{2\pi i j x} \Biggr| \\
&= \sup_{0\le x\le 1}
   \Biggl| X_0 + 2\sum_{j=1}^{n-1} \left(1-\frac{j}{n}\right)
                   X_j \cos(2\pi j x) \Biggr| \\
&= \sup_{0\le x\le 1}\Biggl| \sum_{j=0}^{n-1} a_j X_j \varphi_j(x)
   \Biggr|,
\end{align*}
where we have defined $a_0=1$, $a_j = \sqrt{2}(1-j/n)$ for $j\ge 1$,
$\varphi_0 \equiv 1$, and $\varphi_j(x) = \sqrt{2}\cos(2\pi j x)$ for
$j\ge 1$.  It is easy to verify that $\|a\|_2 > \sqrt{n}/2$ and $\| a
\|_4 < 2 n^{1/4}$. The theorem now follows from Proposition
\ref{T:KT}.
\end{proof}

We remark that by combining Theorem \ref{T:lower} with the proof of
Theorem \ref{T:limsup}, one obtains a nontrivial bound on the left
tail of $\| T_n \|$ under the assumptions of Theorem \ref{T:limsup} and
the additional assumption that $\E X_j^2=1$ for every
$j$. Unfortunately, one cannot deduce an almost sure lower bound of
the form
\[
\liminf_{n\to\infty} \frac{\|T_n\|}{\sqrt{n\log n}} \ge c
\quad \mbox{almost surely}
\]
without more precise control over the constants in Proposition
\ref{T:KT} and the concentration inequalities used in the proof of
Theorem \ref{T:limsup}.

\bigskip
\section{Extensions and additional remarks}\label{S:remarks}


\subsection{Other random matrix ensembles}
For simplicity Theorems \ref{T:upper}--\ref{T:lower} were stated and
proved only for the case of real symmetric Toeplitz matrices. However,
straightforward adaptations of the proofs show that the theorems
hold for other related ensembles of random matrices. These include
nonsymmetric real Toeplitz matrices $\big[X_{j-k}\big]_{j,k\in\Z}$ for
independent random variables $X_j,$ $j\in \Z$, as well as complex
Hermitian or general complex Toeplitz variants. In the complex cases
one should consider matrix entries of the form $X_j = Y_j + i Z_j$,
where $Y_j$ and $Z_j$ are independent and each satisfy the tail or
moment conditions imposed on $X_j$ in the theorems as stated.

Closely related to the case of nonsymmetric random Toeplitz matrices
are random Hankel matrices $H_n = \big[X_{j+k-1}\big]_{1\le j,k\le
n}$, which are constant along skew diagonals. This ensemble was also
mentioned by Bai \cite{Bai}, and was shown to have a universal
limiting spectral distribution in \cite{BDJ}. Independently, Masri and
Tonge \cite{MT} considered a random $r$-linear Hankel form
\[(v_1,\dotsc,v_r)\mapsto \sum_{j_1,\dotsc,j_r=0}^{n} X_{j_1+\dotsb+j_r}
  (v_1)_{j_1} \dotsb (v_r)_{j_r}
\]
in the case $\Prob [X_j = 1]=\Prob[X_j=-1]=1/2$, and showed that the
expected norm of this form is of the order $\sqrt{n^{r-1} \log n}$.
As observed in \cite[Remark 1.2]{BDJ}, $H_n$ has the same singular
values, and so in particular the same spectral norm, as the
(nonsymmetric) Toeplitz matrix obtained by reversing the order of the
rows of $H_n$. Therefore Theorems \ref{T:upper}--\ref{T:lower} apply
to $H_n$ as well. As mentioned in the introduction, the versions of
Theorems \ref{T:upper} and \ref{T:lower} for $H_n$ generalize the
$r=2$ case of the result of \cite{MT} to subgaussian matrix entries
$X_j$.

The methods of this paper can also be used to treat random Toeplitz
matrices with additional restrictions. For example, the theorems apply
to the ensemble of symmetric circulant matrices considered in
\cite[Remark 2]{BM} which is defined as $T_n$ here except for the
restriction that $X_{n-j}=X_j$ for $j=1,\dotsc,n-1$, and the closely
related symmetric palindromic Toeplitz matrices considered in
\cite{MMS}, in which $X_{n-j-1} = X_j$ for $j=0,\dotsc,n-1$. We remark
that \cite{BM,MMS} show that each of these ensembles, properly scaled
and with some additional assumptions, have a limiting spectral
distribution which is normal.


\medskip

\subsection{Weaker hypotheses}

It is unclear how necessary the tail or moment conditions on the $X_j$
are to the conclusions of the theorems. It appears likely (cf.\
\cite{YBK,BoseSen}) that versions of Theorems \ref{T:upper} and
\ref{T:limsup} remain true assuming only the existence of fourth
moments, at least when the $X_j$ are identically distributed.  In
particular it is very likely that the assumptions of Theorem
\ref{T:limsup} can be relaxed considerably. Even within the present
proof, the assumption of a logarithmic Sobolev inequality can be
weakened slightly to that of a quadratic transportation cost
inequality; cf.\ \cite[Chapter 6]{Ledoux2}.

If the $X_j$ have nonzero means then the behavior of $\|T_n\|$ may
change.  Suppose first that the $X_j$ are uniformly subgaussian and
$\E X_j = m\neq 0$ for every $j$. If $J_n$ denotes the $n\times n$
matrix whose entries are all $1$, then \eqref{E:weak-limsup} implies
that
\begin{equation}\label{E:limsup-m}
\limsup_{n\to\infty} \frac{\|T_n - mJ_n\|}{\sqrt{n}\log n} \le c
\quad \mbox{almost surely,}
\end{equation}
where $c$ depends on $m$ and the subgaussian estimate for the
$X_j$. Since $\|J_n\|=n$, \eqref{E:limsup-m} and the triangle
inequality imply a strong law of large numbers:
\begin{equation}\label{E:BoseSen}
\lim_{n\to \infty} \frac{\|T_n\|}{n} = |m| \quad \mbox{almost surely.}
\end{equation}
In \cite{BoseSen}, \eqref{E:BoseSen} was proved using estimates from
\cite{BDJ} under the assumption that the $X_j$ are identically
distributed and have finite variance. We emphasize again
that while the methods of this paper require stronger tail
conditions, we never assume the $X_j$ to be identically distributed.

More generally, the behavior of $\|T_n\|$ depends on the rate of
growth of the spectral norms of the deterministic Toeplitz matrices
$\E T_n$.  The same argument as above shows that
\[
\lim_{n\to\infty} \frac{\|T_n\|}{\| \E T_n \|} = 1 \quad 
\mbox{almost surely}
\]
if the random variables $(X_j-\E X_j)$ are uniformly subgaussian and
$\lim_{n\to\infty} \frac{\sqrt{n}\log n}{\| \E T_n \|} = 0$.
On the other hand, if $\| \E T_n\| = o(\sqrt{n\log n})$ then the
conclusion of Theorem \ref{T:upper} holds.

\subsection{Random trigonometric polynomials}

The supremum of the random trigonometric polynomial
\[
Z_x = \sum_{j=1}^n X_j \cos(2\pi j x),
\]
has been well-studied in the special case $\Prob[X_j=1] =
\Prob[X_j=-1] = 1/2$, in work dating back to Salem and Zygmund
\cite{SZ}. Observe that $Z_x$ is essentially equivalent to the process
$Y_x$ defined in the proof of Theorem \ref{T:upper}, and is also
closely related to the random process considered in the proof of
Theorem \ref{T:lower}. Hal\'asz \cite{Halasz} proved in particular that
\[
\lim_{n\to \infty}
\frac{\sup_{0\le x\le 1} |Z_x|}{\sqrt{n\log n}}
= 1 \quad \mbox{almost surely}.
\]
From this it follows that when $\Prob[X_j=1] = \Prob[X_j=-1] = 1/2$
for every $j$, the conclusion of Theorem \ref{T:limsup} holds with
$c_2 = 2$. Numerical experiments suggest, however, that the optimal
value of $c_2$ is $1$ in this case, and more generally when the $X_j$
are i.i.d.\ with mean $0$ and variance $1$.

Conversely, adaptations of the proofs in this paper yield less
numerically precise bounds for the supremum of $Z_x$ under the same
weaker assumptions on the $X_j$ in the statements of the theorems. We
remark that the techniques used to prove the results of
\cite{KT1,KT2,MT} cited above (and hence indirectly also Theorem
\ref{T:lower}) were adapted from the work of Salem and Zygmund in
\cite{SZ}.


\bibliographystyle{plain} \bibliography{snrtm}

\begin{thebibliography}{10}

\bibitem{Bai}
Z.~D. Bai.
\newblock Methodologies in spectral analysis of large-dimensional random
  matrices, a review.
\newblock {\em Statist. Sinica}, 9(3):611--677, 1999.

\bibitem{BM}
A.~Bose and J.~Mitra.
\newblock Limiting spectral distribution of a special circulant.
\newblock {\em Statist. Probab. Lett.}, 60(1):111--120, 2002.

\bibitem{BoseSen}
A.~Bose and A.~Sen.
\newblock Spectral norm of random large dimensional noncentral {T}oeplitz and
  {H}ankel matrices.
\newblock {\em Electron. Comm. Probab.}, 12:29--35, 2007.

\bibitem{BS}
A.~B{\"o}ttcher and B.~Silbermann.
\newblock {\em Introduction to Large Truncated {T}oeplitz Matrices}.
\newblock Universitext. Springer-Verlag, New York, 1999.

\bibitem{BDJ}
W.~Bryc, A.~Dembo, and T.~Jiang.
\newblock Spectral measure of large random {H}ankel, {M}arkov and {T}oeplitz
  matrices.
\newblock {\em Ann. Probab.}, 34(1):1--38, 2006.

\bibitem{Dudley}
R.~M. Dudley.
\newblock The sizes of compact subsets of {H}ilbert space and continuity of
  {G}aussian processes.
\newblock {\em J. Funct. Anal.}, 1:290--330, 1967.

\bibitem{Halasz}
G.~Hal{\'a}sz.
\newblock On a result of {S}alem and {Z}ygmund concerning random polynomials.
\newblock {\em Studia Sci. Math. Hungar.}, 8:369--377, 1973.

\bibitem{HM}
C.~Hammond and S.~J. Miller.
\newblock Distribution of eigenvalues for the ensemble of real symmetric
  {T}oeplitz matrices.
\newblock {\em J. Theoret. Probab.}, 18(3):537--566, 2005.

\bibitem{KT1}
B.~Kashin and L.~Tzafriri.
\newblock Lower estimates for the supremum of some random processes.
\newblock {\em East J. Approx.}, 1(1):125--139, 1995.

\bibitem{KT2}
B.~Kashin and L.~Tzafriri.
\newblock Lower estimates for the supremum of some random processes, {II}.
\newblock {\em East J. Approx.}, 1(3):373--377, 1995.

\bibitem{Ledoux1}
M.~Ledoux.
\newblock Concentration of measure and logarithmic {S}obolev inequalities.
\newblock In {\em S\'eminaire de Probabilit\'es, XXXIII}, volume 1709 of {\em
  Lecture Notes in Math.}, pages 120--216. Springer, Berlin, 1999.

\bibitem{Ledoux2}
M.~Ledoux.
\newblock {\em The Concentration of Measure Phenomenon}, volume~89 of {\em
  Mathematical Surveys and Monographs}.
\newblock American Mathematical Society, Providence, RI, 2001.

\bibitem{LPRT}
A.~E. Litvak, A.~Pajor, M.~Rudelson, and N.~Tomczak-Jaegermann.
\newblock Smallest singular value of random matrices and geometry of random
  polytopes.
\newblock {\em Adv. Math.}, 195(2):491--523, 2005.

\bibitem{MT}
I.~Masri and A.~Tonge.
\newblock Norm estimates for random multilinear {H}ankel forms.
\newblock {\em Linear Algebra Appl.}, 402:255--262, 2005.

\bibitem{MMS}
A.~Massey, S.~J. Miller, and J.~Sinsheimer.
\newblock Distribution of eigenvalues of real symmetric palindromic {T}oeplitz
  matrices and circulant matrices.
\newblock {\em J. Theoret. Probab.}
\newblock To appear. Preprint available at {\tt
  http://arxiv.org/math.PR/0512146}.

\bibitem{Rudelson}
M.~Rudelson.
\newblock Probabilistic and combinatorial methods in analysis.
\newblock Lecture notes from an NSF-CBMS Regional Research Conference at Kent
  State University, 2006.

\bibitem{SZ}
R.~Salem and A.~Zygmund.
\newblock Some properties of trigonometric series whose terms have random
  signs.
\newblock {\em Acta Math.}, 91:245--301, 1954.

\bibitem{Talagrand2}
M.~Talagrand.
\newblock Majorizing measures: the generic chaining.
\newblock {\em Ann. Probab.}, 24(3):1049--1103, 1996.

\bibitem{YBK}
Y.~Q. Yin, Z.~D. Bai, and P.~R. Krishnaiah.
\newblock On the limit of the largest eigenvalue of the large-dimensional
  sample covariance matrix.
\newblock {\em Probab. Theory Related Fields}, 78(4):509--521, 1988.

\end{thebibliography}

\end{document}